\documentclass[a4paper,12pt]{article}

\usepackage{amsfonts}
\usepackage{amscd,color}
\usepackage{amsmath,amsfonts,amssymb,amscd}
\usepackage{indentfirst,graphicx,epsfig}
\usepackage{tikz}
\usepackage{graphicx}
\input{epsf}
\usepackage{epstopdf}
\usepackage{caption}

\newtheorem{thm}{Theorem}[section]

\newtheorem{lem}[thm]{Lemma}

\newcommand{\ml}{l\kern-0.55mm\char39\kern-0.3mm}

\def\qed{\hfill \nopagebreak\rule{5pt}{8pt}}

\title{\textbf{Hardness results for rainbow disconnection of graphs\footnote{Supported by NSFC No.11871034, 11531011 and NSFQH No.2017-ZJ-790.}}}
\author{\small \ Zhong Huang, \ Xueliang Li \\
\small $^1$Center for Combinatorics and LPMC \\
\small Nankai University, Tianjin 300071, China\\
\small $^2$School of Mathematics and Statistics\\
\small Qinghai Normal University\\
\small Xining, Qinghai 810008, China\\
{\small 2120150001@mail.nankai.edu.cn, lxl@nankai.edu.cn}
}
\date{}
\begin{document}
\maketitle
\begin{abstract}
Let $G$ be a nontrivial connected, edge-colored graph. An edge-cut $S$ of $G$
is called a rainbow cut if no two edges in $S$ are colored with a same color.
An edge-coloring of $G$ is a rainbow disconnection coloring if for every two distinct
vertices $s$ and $t$ of $G$, there exists a rainbow cut $S$ in $G$ such that $s$ and $t$
belong to different components of $G\setminus S$. For a connected graph $G$, the {\it rainbow disconnection
number} of $G$, denoted by $rd(G)$, is defined as the smallest number of colors such that $G$ has a rainbow disconnection coloring by using this number of colors. In this paper, we show that
for a connected graph $G$, computing $rd(G)$ is NP-hard. In particular, it is already NP-complete to
decide if $rd(G)=3$ for a connected cubic graph. Moreover, we prove that for a given edge-colored
(with an unbounded number of colors) connected graph $G$ it is NP-complete to decide
whether $G$ is rainbow disconnected.    \\[2mm]
\end{abstract}

\section{Introduction}

All graphs in this paper are simple, finite and undirected. We
follow \cite{BM} for graph theoretical notation and terminology not
described here. Let $G$ be a graph. We use $V(G), E(G), n(G), m(G)$, $\delta(G)$
and $\Delta(G)$ to denote the vertex-set, edge-set, number of
vertices, number of edges, minimum degree and maximum degree of $G$, respectively.
Let $c: E(G) \rightarrow [k] = \{1, 2, ..., k\}, k \in N$ be an edge-coloring of $G$,
where adjacent edges may be colored with a same color. When adjacent edges of $G$ receive
different colors under $c$, the edge-coloring $c$ is called {\it proper}. The {\it chromatic index} of $G$,
denoted by $ \chi^\prime(G)$, is the minimum number of colors needed in a proper coloring of
$G$. By a famous theorem of Vizing \cite{VG} we have
\[\Delta(G) \leq\chi^\prime(G) \leq \Delta(G) + 1\]
for every nonempty graph $G$. And, if $\chi^\prime(G)= \Delta(G)$, then $G$ is of Class $1$;
if $\chi^\prime(G) =\Delta(G) + 1$, then $G$ is of Class $2$.

A path of an edge-colored graph $G$ is said to be a
\emph{rainbow path} if no two edges on the path have the same color.
The graph $G$ is called \emph {rainbow connected} if every pair
of distinct vertices of $G$ is connected by a rainbow path in $G$.
An edge-coloring of a connected graph is a \emph{rainbow connection
coloring} if it makes the graph rainbow connected. This concept of
rainbow connection of graphs was introduced by Chartrand et
al.~\cite{CJMZ} in 2008. For a connected graph $G$, the \emph{rainbow connection number}
$rc(G)$ of $G$ is defined as the smallest number of colors
that are needed in order to make $G$ rainbow connected. The reader
who are interested in this topic can see \cite{LSS, LS} for a survey.

An {\it edge-cut} of a nontrivial connected graph $G$ is a set $S$ of edges of $G$ such that $G \setminus S$ is disconnected.
The minimum number of edges in an edge-cut of $G$ is defined as the {\it edge-connectivity} $\lambda(G)$ of $G$.
We then have the well-known inequality $\lambda(G)\leq \delta(G)$. For two distinct vertices $s$ and $t$ of $G$,
let $\lambda(s, t)$ denote the minimum number of edges in an edge-cut $S$ of $G$ such that $s$ and $t$ lie in different components
of $G \setminus S$. The so-called {\it upper edge-connectivity} $\lambda^+(G)$ of $G$ is defined by
\[\lambda^+(G) = max\{\lambda(s, t) : s, t \in V(G)\}.\]
$\lambda(G)$ denotes the global minimum
edge-connectivity of a graph, while $\lambda^+(G)$ denotes the local maximum edge-connectivity of a graph.

An edge-cut $S$ of an edge-colored graph $G$ is called a {\it rainbow cut} if no two edges in $S$ are colored
with a some color. A rainbow cut $S$ is said to {\it separate} two vertices $s$ and $t$ if $s$ and $t$ belong to
different components of $G\setminus S$. Such a rainbow cut is called a {\it $s-t$ rainbow cut}. An edge-colored graph
$G$ is called a {\it rainbow disconnected} if for every two vertices $s$ and $t$ of $G$, there exists an $s-t$ rainbow cut in $G$.
In this case, the edge-coloring $c$ is called a {\it rainbow disconnection coloring} of $G$.
Similarly, we define the {\it rainbow disconnection number} of $G$, denoted by $rd(G)$, as the smallest number of
colors such that $G$ has a rainbow disconnected coloring by using this number of colors. A rainbow disconnection coloring with $rd(G)$
colors is called an {\it $rd$-coloring} of $G$. This concept of rainbow disconnection of graphs was introduced by Chartrand
et al.~\cite{GST} very recently in 2018.

In this paper, we show that for a connected graph $G$ computing $rd(G)$ is NP-hard.
In particular, it is already NP-complete to decide if $rd(G)=3$ for a connected cubic graph $G$.
Moreover, we show that for a given edge-colored (with an unbounded number of colors) connected graph $G$
it is NP-complete to decide whether $G$ is rainbow disconnected.

\section{Hardness results}

At the very beginning, we state some fundamental results on the
rainbow disconnection of graphs, which will be used in the sequel.

\begin{lem}{\upshape \cite{GST}}\label{lem2-1}
If $G$ is a nontrivial connected graph, then
\[\lambda(G)\leq \lambda ^+
(G) \leq rd(G) \leq \chi^\prime(G)\leq \Delta(G) + 1.\]
\end{lem}

Next result is due to Holyer \cite{HI}, which is on the complexity of the chromatic index of a cubic graph.

\begin{lem}{\upshape \cite{HI}}\label{lem2-2}
It is NP-complete to determine whether the chromatic index of a cubic graph is 3 or 4.
\end{lem}

At first we should show that our problem is in NP for a fixed integer $k$.

\begin{lem}\label{lem2-3}
For a fixed positive integer $k$, given a $k$-edge-colored graph
$G$, deciding whether $G$ is rainbow disconnected is in $P$.
\end{lem}

\begin{pf}
Let $n, m$ be the number of vertices and edges of $G$ respectively. Let $s, t$ be two vertices in $G$. Since $G$ is $k$-edge-colored, any rainbow cut set $S$ contains at most $k$ edges, and so, we have no more than $C_m^k$ choices of $S$. Given a set $S$ of edges, it is easy to check whether $s$ and $t$ lie in different components of $G\setminus S$. And there are at most $C_n^2$ pairs of vertices in $G$. Then, we can deduce that deciding whether $G$ is rainbow disconnected can be checked in polynomial-time.
\end{pf}

Next lemma is crucial for proof of our main result.

\begin{lem}\label{lem2-4}
Let $G$ be a $3$-edge-connected cubic graph. Then $\chi^\prime (G)=3$ if and only if $rd(G)=3$.
\end{lem}

\begin{pf}
Assume that $\chi^\prime (G)=3$, and let us show that $rd(G)=3$. Noticing that $G$ is $3$-edge-connected, we have that $rd(G)\geq 3$. Since $rd(G)\leq \chi^\prime (G)$ by Lemma \ref{lem2-1}, we then have $rd(G)=3$.

Assume that $rd(G)=3$. Let $S=\{u_1v_1,u_2v_2,u_3v_3\}$ be a rainbow $3$-edge cut of $G$, and $G\setminus S$ has two non-trivial component (which means the component is not a singleton) $C_1$ and $C_2$. If all the edges of $S$ share a common vertex, then one of $C_1$ of $C_2$ is a singleton, a contradiction. If two edges of $S$ are adjacent, say $u_1=u_2$, let $e$ be the third edge which is adjacent to $u_1$, then $S^\prime=\{e,u_3v_3\}$ is a $2$-edge cut of $G$, a contradiction. If non edges of $S$ are adjacent, then we employ a new vertex $x_1$ which is adjacent to $u_1,u_2,u_3$ in $C_1$, and $u_ix_1$ receive the same color as $u_iv_i$ for $i=1,2,3$. Similarly, we employ a new vertex $x_2$ which is adjacent to $v_1,v_2,v_3$ in $C_2$, and $x_2v_i$ receive the same color as $u_iv_i$ for $i=1,2,3$. Now we get $3$-edge connected cubic graphs $C^\prime_1$ and $C^\prime_2$. Repeat the operation on $C^\prime_1$ and $C^\prime_2$, and finally we get a graph sequence $T=\{T_1,T_2 \cdots T_r\}$. Let $s, t$ be two vertices in $T_j\in T$, and $S_{s,t}$ be the rainbow cut, we then have that $|S_{s,t}|=3$ and the three edges of $S_{s,t}$ are incident with one of $s, t$. Next, we deduce that every vertex of $T_j$ is incident with three rainbow colored edges except for one vertex, say $s_0$. Let the number of edges with color $i$ incident with $s_0$ be $k_i$, and $T_{12}$ be the subgraph induced by the set of edges of $T_j$ which are colored with colors $1$ or $2$. Let $s_{12}\neq s_0$ be a vertex of $T_{12}$, then $d(s_{12})=2$. Since the degree sum of $T_{12}$ is an even number, we have $k_1+k_2=0 \pmod 2$, which gives $k_1 = k_2 \pmod 2$. Similarly, $k_2 = k_3 \pmod 2$. So, we have that $k_1=k_2=k_3=1$ and $s_0$ is incident with three rainbow colored edges. So, $T_j$ is properly colored for $1\leq j\leq r$. Now we consider the original graph $G$, the coloring satisfies that $rd(G)=3$ is also a proper edge-coloring. Then, $\chi^\prime (G)=3$.\qed
\end{pf}

In the proof of the following result, we will use the graph $G_\phi$ which contains some copies of Figures 1, 2 and 3 employed in Holyer's paper \cite{HI}.

\begin{figure}[!htb]
\centering
\includegraphics[width=0.8\textwidth]{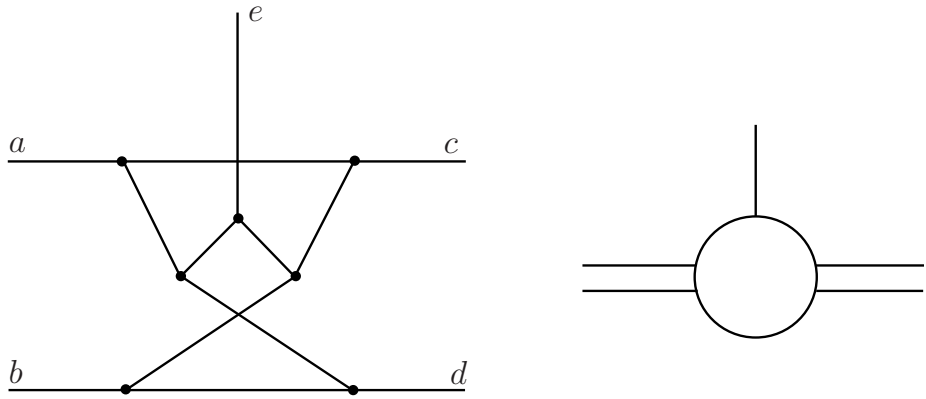}
\caption{The inverting component and its symbolic representation.}
\end{figure}

\begin{figure}[!htb]
\centering
\includegraphics[width=0.8\textwidth]{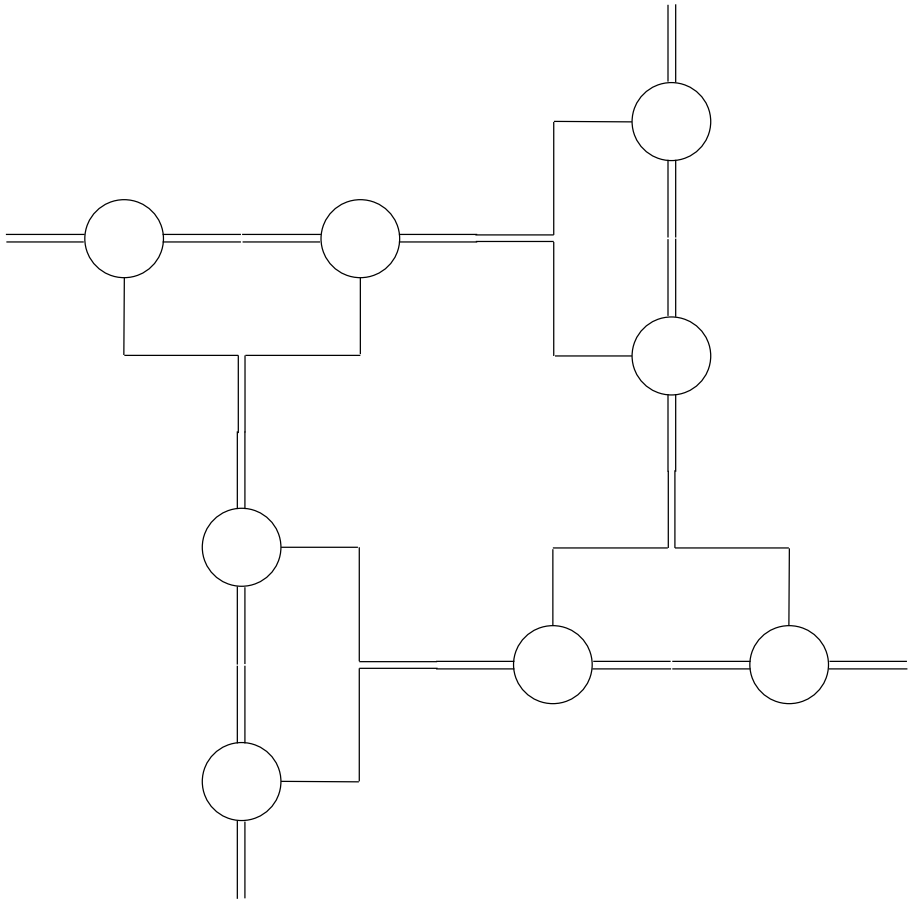}
\caption{The variable-setting component made from 8 inverting components and having 4 output pairs of edges. More generally, it is made from $2n$ inverting components and has $n$ output pairs $(n\geq 2)$.}
\end{figure}

\begin{figure}[!htb]
\centering
\includegraphics[width=0.8\textwidth]{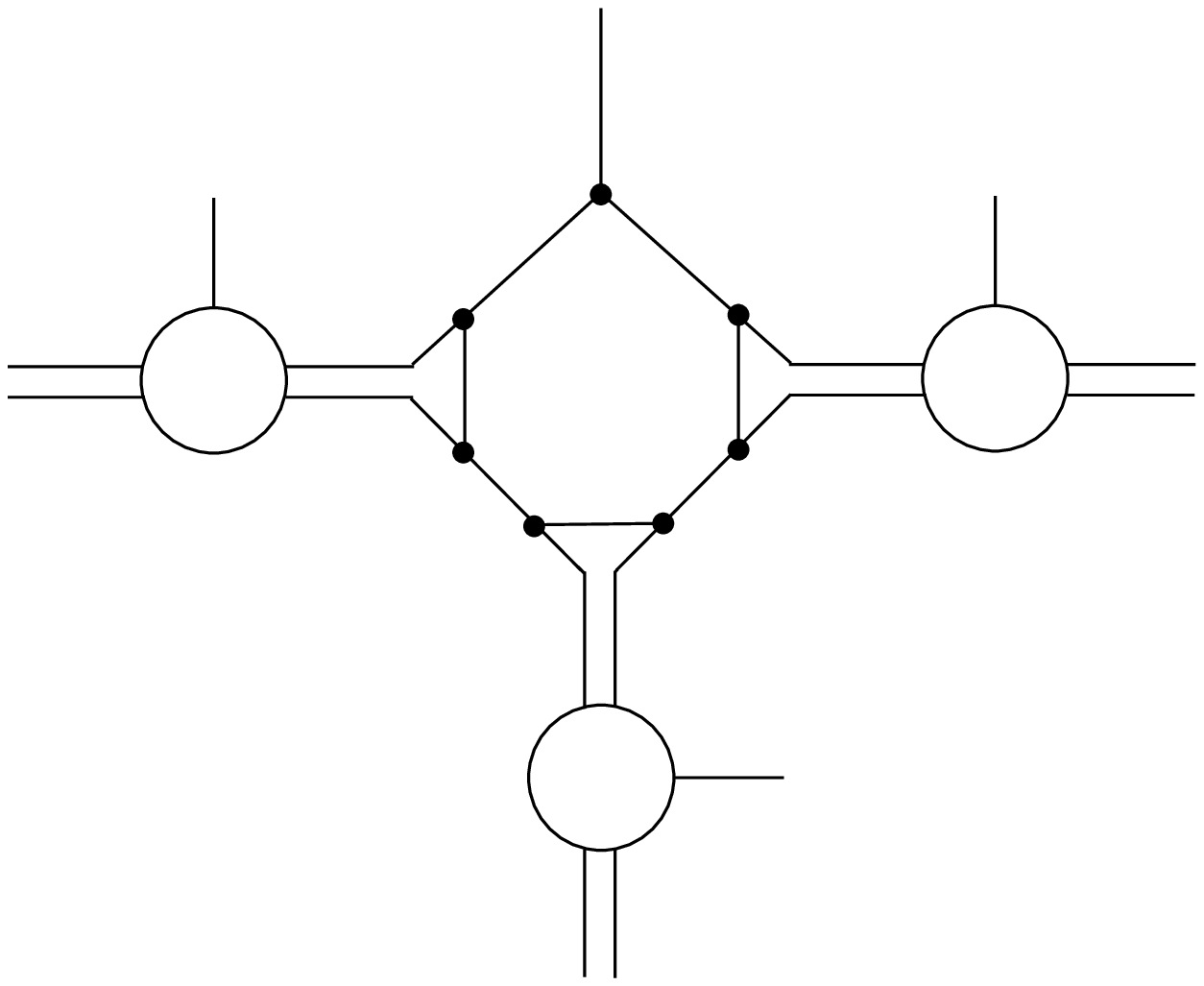}
\caption{The satisfaction testing component.}
\end{figure}

\begin{thm}\label{th2-1}
It is NP-complete to determine whether the rainbow disconnection number of a cubic graph is $3$ or $4$.
\end{thm}

\begin{pf}
Clearly, the problem is in NP by Lemma \ref{lem2-3}. We prove that it is NP-complete by reducing $3$-SAT to it. Given a $3$CNF formula $\phi=\wedge_{i=1}^m C_i$ over $n$ variables $x_1, x_2, \cdots, x_n$, we use the cubic graph $G_\phi$ that was used by Holyer in Lemma \ref{lem2-2} such that $rd(G_\phi)=3$ if and only if $\phi$ is satisfiable.

Noticing that $G(\phi)$ is $3$-edge-connected, we then can verify that $rd(G_\phi)=3$ if and only if $\phi$ is satisfiable by Lemma \ref{lem2-4}.\qed
\end{pf}

Deciding whether a $k$-edge-colored graph
$G$, where $k$ is a constant, is rainbow disconnected is in $P$. However, it is NP-complete to decide whether a given edge-colored (with an unbounded number of colors) graph is rainbow disconnected.
The proof of the following result uses a similar technique of \cite{CJMZ}.

\begin{thm}\label{th2-2}
Given an edge-colored graph $G$ and two vertices $s,t$ of $G$, deciding whether there is a rainbow cut between $s$ and $t$ is NP-complete.
\end{thm}

\begin{pf}
Clearly, the problem is in NP. We prove that it is NP-complete by reducing $3$-SAT to it. Given a $3$CNF formula $\phi=\wedge_{i=1}^m c_i$ over $n$ variables $x_1,x_2,\cdots,x_n$, we construct a graph $G_{\phi}$ with two special vertices $s,t$ and a coloring $c: E(G_{\phi})\rightarrow [E(G_{\phi})]$ such that there is a rainbow cut between $s$ and $t$ in $G_{\phi}$ if and only if $\phi$ is satisfiable.

We define $G_{\phi}$ as follows:
\[
V(G_{\phi})=\{c_i,c_i^1,c_i^2,c_i^3: i\in [m]\}\cup \{x_i^0,x_i^1: i\in[n]\}\cup \{s,t\}
\]
\[
\begin{aligned}
E(G_{\phi})=&\Big\{\{x_j^0,c_i\},\{c_i,c_i^k\},\{c_i^k,x_j^1\}:\\
&\text{If $x_j\in c_i$ corresponds to $c_i^k$ is positive in $\phi$},  k\in[1,2,3] \Big\}\\
&\cup \Big\{\{x_j^1,c_i\},\{c_i,c_i^k\},\{c_i^k,x_j^0\}:\\
&\text{If $x_j\in c_i$ corresponds to $c_i^k$ is negative in $\phi$}, k\in[1,2,3] \Big\} \\
&\cup \Big\{\{s,x_i^0\},\{s,x_i^1\}: i\in[n]\}\Big\}\cup \Big\{ \{s,t\} \Big\}\cup\Big\{ \{t,c_i\}:i\in[m] \Big\}
\end{aligned}
\]
The coloring $c$ is defined as follows:

\begin{itemize}
\item the edges $\Big\{\{s,t\},\{t,c_i\}:i\in[m]\Big\}$ are colored with a special color $r_0$;
\item the edges $\Big\{\{s,x_i^0\},\{s,x_i^1\}: i\in[n]\}\Big\}$ are colored with a special color $r_i$, $i\in[n]$;
\item the edge $\{x_j^0,c_i\}$ or $\{x_j^1,c_i\}$ is colored with a special color $r_i^k$ when $x_j\in c_i$ is positive or negative in $\phi$ respectively, $i\in [m], k\in[1,2,3]$;
\item the edge $\{c_i^k,x_j^1\}$ or $\{c_i^k,x_j^0\}$ is colored with a special color $r_i^4$ when $x_j\in c_i$ is positive or negative in $\phi$ respectively, $i\in [m], k\in[1,2,3]$;
\item the edges $\Big\{\{c_i, c_i^j\}:j\in[1,2,3]\Big\}$ are colored with a special color $r_i^5$, $i\in [m]$.
\end{itemize}
Now we can verify that there is a rainbow cut between $s$ and $t$ in $G_{\phi}$ if and only if $\phi$ is satisfiable.

Assume that there is a rainbow cut $S$ between $s$ and $t$ in $G_\phi$ under $c$, and let us show that $\phi$ is satisfiable. At first, we consider the color $r_0$. Since $s$ and $t$ are adjacent in $G(\phi)$, then the edge $\{s,t\}$ is in $S$. Clearly, $S$ separates $s$ and the set $\{t\}\cup \{c_i:i\in[m]\}$. Next, the color $r_i$ appears twice in $G(\phi)$. Without loss of generality, we can assume that there is exact one of $\{s,x_i^0\}$ and $\{s,x_i^1\}$ in $S$, which corresponds to the value of variable $x_i$, $i\in[n]$. At last, we consider the colors $r_i^4$ and $r_i^5$, $i\in[m]$. There are at most two edges that have color $r_i^4$ or $r_i^5$ in $S$, which means that $c_i$ (a clause in $\phi$) is satisfiable, $i\in[m]$. As a result, $\phi$ is satisfiable.

Assume that $\phi$ is satisfiable, and let us construct a rainbow cut $S$ between $s$ and $t$ in $G_\phi$ under $c$. At first, the edges $\{s,x_i^{|x_i|}\},i\in[n]$ and \{s,t\} are in $S$. If the vertex $x_j^{|x_j|}$ is adjacent to $c_i$, then we choose one edge colored with $r_i^4$ or $r_i^5$ corresponding to variable $x_j$. If the vertex $x_j^{|x_j|}$ is adjacent to $c_i^k$, then we choose the edge
$\{c_i,x_j^{|\overline x_j|}\}$ that is colored with $r_i^k$ and corresponds to variable $x_j$. Notice that $\phi$ is satisfiable and no more than two edges colored with $r_i^4$ or $r_i^5$ are chosen. Add these chosen edges to $S$, and now $S$ is a rainbow cut between $s$ and $t$ in $G_\phi$ under $c$.\qed
\end{pf}

\end{document}